\documentclass[10pt, a4paper]{amsart}
\usepackage{amsthm}
\usepackage{graphicx}
\theoremstyle{}
{\theoremstyle{definition}
\newtheorem{dfn}{Definition}[section]}
\newtheorem{prop}[dfn]{Proposition}

\newtheorem{thm}[dfn]{Theorem}
{\theoremstyle{definition}
\newtheorem{rem}[dfn]{Remark}}

{\theoremstyle{definition}
}

\usepackage[top=30truemm,bottom=30truemm,left=30truemm,right=30truemm]{geometry}
\usepackage{amsmath,amssymb,amscd,bm,graphicx,tikz-cd}
\usepackage[all]{xy}
\usepackage{hyperref}

\newcommand{\bC}{\mathbb{C}}
\newcommand{\bF}{\mathbb{F}}

\newcommand{\bH}{\mathbb{H}}

\newcommand{\bZ}{\mathbb{Z}}

\begin{document}

\title[]{Twisted Alexander Polynomials of $(-2,3,2n+1)$-Pretzel Knots}

\author[A.~Aso]{Airi Aso}
\email{aso-airi@ed.tmu.ac.jp}
\address{Department of mathematics and information sciences, Tokyo Metropolitan University, 1-1 Minamiohsawa, Hachioji-shi, Tokyo, 192-0397, Japan}

\begin{abstract} 
We calculate the twisted Alexander polynomials of $(-2,3,2n+1)$-pretzel knots associated to their holonomy representations. As a corollary, we obtain new supporting evidences of Dunfield, Friedl and Jackson's conjecture, that is, the twisted Alexander polynomials of hyperbolic knots associated to their holonomy representations determine the genus and fiberedness of the knots.
\end{abstract}

\maketitle{}

\section{Introduction}
\subsection{Motivations}
The notion of twisted Alexander polynomials was introduced by Wada [W] and Lin [L] independently in 1990s.
The definition of Lin is for knots in $S^3$ and the definition of Wada is for finitely presented groups. 
The twisted Alexander polynomial is a generalization of the Alexander polynomial, and it is defined for the pair of a group and its representations.
By Kitano and Morifuji [KM], it is known that Wada's twisted Alexander polynomials of the knot groups for any nonabelian representations into $SL_2(\bF)$ over a field $\bF$ are polynomials. 
As a corollary of the claim, they also showed that if $K$ is a fibered knot of genus $g$,
then its twisted Alexander polynomials are monic polynomials of degree $4g-2$ for any nonabelian $SL_2(\bF)$-representations. The converse does not hold, in other words, there exist examples of nonfibered knot which has a $SL_2(\bC)$-representation such that the twisted Alexander polynomial of the representation is monic (see [GoMo]).  

If $K$ is hyperbolic, i.e. the complement $S^3 \setminus K$ admits a complete hyperbolic metric of finite volume,
the most important representation is its holonomy representation into $SL_2(\bC)$ which is a lift of the representation into the group of orientation-preserving isometries of the hyperbolic 3-space $\bH^3$.
Dunfield, Friedl and Jackson [DFJ] conjectured that the twisted Alexander polynomials of hyperbolic knots associated to their holonomy representations determine the genus and fiberedness of the knots.
In fact, they computed the twisted Alexander polynomials of all hyperbolic knots of 15 or fewer crossings associated to their holonomy representations, and the conjecture is verified for these hyperbolic knots. 
Recently, the twisted Alexander polynomials of some infinite families of knots, twist knots and genus one two-bridge knots associated to their holonomy representations, are computed by Morifuji [Mo1] and Tran [T1], and genus one two-bridge knots associated to the adjoint representations of their holonomy representations is also computed by Tran [T2]. 
These examples are also the supporting evidences of the conjecture.

In this paper, we compute the twisted Alexander polynomials of $(-2,3,2n+1)$-pretzel knots $K_n$ depicted in Figure 1 associated to their holonomy representations given in the following section.
As a corollary, we obtain new supporting evidences of Dunfield, Friedl and Jackson's conjecture, i.e. the twisted Alexander polynomials of $K_n$ are monic polynomials of degree $4n+6$. 
We can observe that $K_n$ is fibered for any non-negative integers $n$ and the genus of $K_n$ is $n+2$ because the Seifert surface $S_n$ is a Murasugi sum of a Seifert surface of a torus knot and a Seifert surface of a Hopf link. Since these Seifert surfaces are fibered surfaces, $S_n$ is also a fibered surface (see [HM, M, O] for more details).

On the other hand, $(-2,3,2n+1)$-pretzel knot is an infinite family of knots which contains the Fintushel-Stern knot i.e. $(-2,3,7)$-pretzel knot. It plays an important role in studying of exceptional surgeries of knots [Ma].
In fact, the A-polynomials of $(-2,3,2n+1)$-pretzel knot are computed by Tamura-Yokota [TY] and Garoufalidis-Mattman [GaMa].

\begin{figure}[h]
  \begin{center}
\includegraphics[clip,width=7cm]{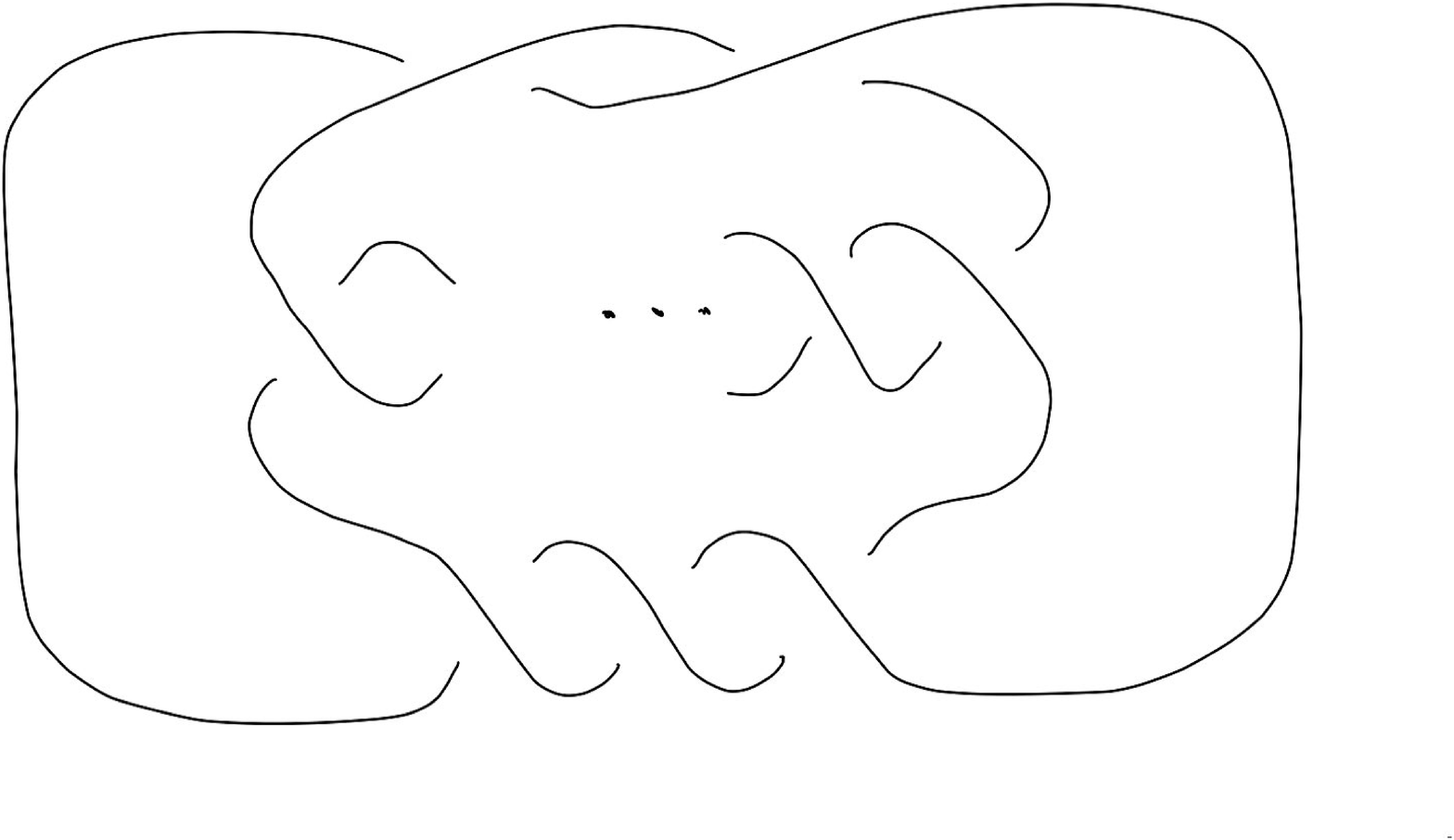}\caption{$(-2,3,2n+1)$-pretzel knot}
  \end{center}
\end{figure}

\subsection{Definition of twisted Alexander polynomials}
In this paper, we use the following definition due to Wada.

\begin{dfn}
Let $G(K)=\pi_1(S^3 \setminus K)$ be the knot group of a knot $K$ presented by
\[
G(K) = \langle x_1, \cdots ,x_n \ \vline \  r_1, \cdots , r_{n-1}  \rangle.
\]
Let $\Gamma$ denote the free group generated by $x_1, \cdots ,x_n $ 
and $\phi: \bZ \Gamma \to \bZ G(K)$ the natural ring homomorphism.
Let $\rho: G(K) \to GL_d(\bF)$ be a $d$-dimensional  linear representation of $G(K)$
and $\Phi : \bZ \Gamma \to M_d(\bF[t,t^{-1}])$ the ring homomorphism defind by
\[
\Phi =(\tilde{\rho} \otimes \tilde{\alpha}) \circ  \phi,  
\]
where $\tilde{\alpha}:\bZ G(K) \to \bZ \langle t,t^{-1} \rangle$ and $\tilde{\rho}$ are respective ring homomorphisms induced by the abelianization $\alpha: G(K) \to  \langle t \rangle$ and $\rho$.
We put
\[
A_{i,j} = \Phi \left( \frac{\partial r_i}{\partial x_j}\right),
\]
where $\displaystyle \frac{\partial}{\partial x_j}$ denotes the Fox derivative (or free derivative) with respect to $x_j$, that is, a map $\bZ \Gamma \to \bZ \Gamma$ satisfying the conditions
\begin{eqnarray*}
\displaystyle \frac{\partial}{\partial x_j} x_i = \delta_{ij} \ \ \ \mbox{and} \ \ \
\displaystyle \frac{\partial}{\partial x_j}g g'= \displaystyle \frac{\partial}{\partial x_j} g+\displaystyle \frac{\partial} {\partial x_j}g',
\end{eqnarray*}
where $\delta_{ij}$ denotes the Kronecker symbol and $g,g' \in \Gamma$.
Then, the twisted Alexander polynomial of $K$ is defined by
\[
\Delta_{K,\rho} =\displaystyle \frac{\det A_{\rho,k}}{\det \Phi(x_k-1)},
\]
where $A_{\rho,k}$ is the $2(n-1)\times 2(n-1)$ matrix obtained from $A_{\rho}=(A_{i,j})$ by removing the $k$-th column, i.e.
\[
A_{\rho,k}=
\left(
\begin{array}{cccccc}
 A_{1,1}& \cdots &  A_{1,k-1} & A_{1,k+1}  & \cdots & A_{1,n}\\
  \vdots & & \vdots & \vdots & &\vdots\\
 A_{n-1,1}& \cdots &  A_{n-1,k-1} & A_{n-1,k+1}  & \cdots & A_{n-1,n}
 \end{array}
\right).
\]
\end{dfn}

\vspace{3mm}
\hspace{-5.3mm}
{\it Acknowledgement}:
The author would like to thank professor Yoshiyuki Yokota for supervising and giving helpful comments. She also would like to thank professor Teruhiko Soma and professor Manabu Akaho for giving valuable comments.

\section{Holonomy representations}
In this section, we give a presentation of knot group $G(K_n)$ and its holonomy representation $\rho_m : G(K_n) \to SL_2(\bC)$, where $m$ represents the eigenvalue of the meridian of $K_n$.

Let $L$ be the link depicted in Figure 2 and $E=S^3 \setminus L$. Then, the Wirtinger presentation (see [CF]) of $\pi_1(E)$ is given by
\[
\langle a,b,x \ \vline \ \{a x b a (x b)^{-1} \}^{-1} x = x b \{a x b a (x b)^{-1}\}^{-1} (a x b)^{-1} x b, \  [x, a x b a (x b)^{-1} ] = 1 \rangle,
\]
where $a,b$ and $x$ is Wirtinger generators assigned to the corresponding pass depicted in Figure 2.
Note that $E_n:=S^3 \setminus K_n$ is obtained from $L$ by $(-\frac{1}{n})$-surgery along the trivial component, that is, removing the tubular neighborhood of the trivial component and re-gluing the solid torus again. Therefore, by the van Kampen theorem, we have
 \begin{eqnarray*}
\pi_1 (E_n) \!\!&\!\! = \!\!&\!\! \langle a,b,x \ \vline \ \{ a x b a (x b)^{-1} \}^{-1} x = x b \{a x b a (x b)^{-1}\}^{-1} (a x b)^{-1} x b, \ 
x =  \{a x b a (x b)^{-1} \}^{n} \rangle.
 \end{eqnarray*}

\begin{figure}[htbp]
  \begin{center}
\includegraphics[clip,width=6.1cm]{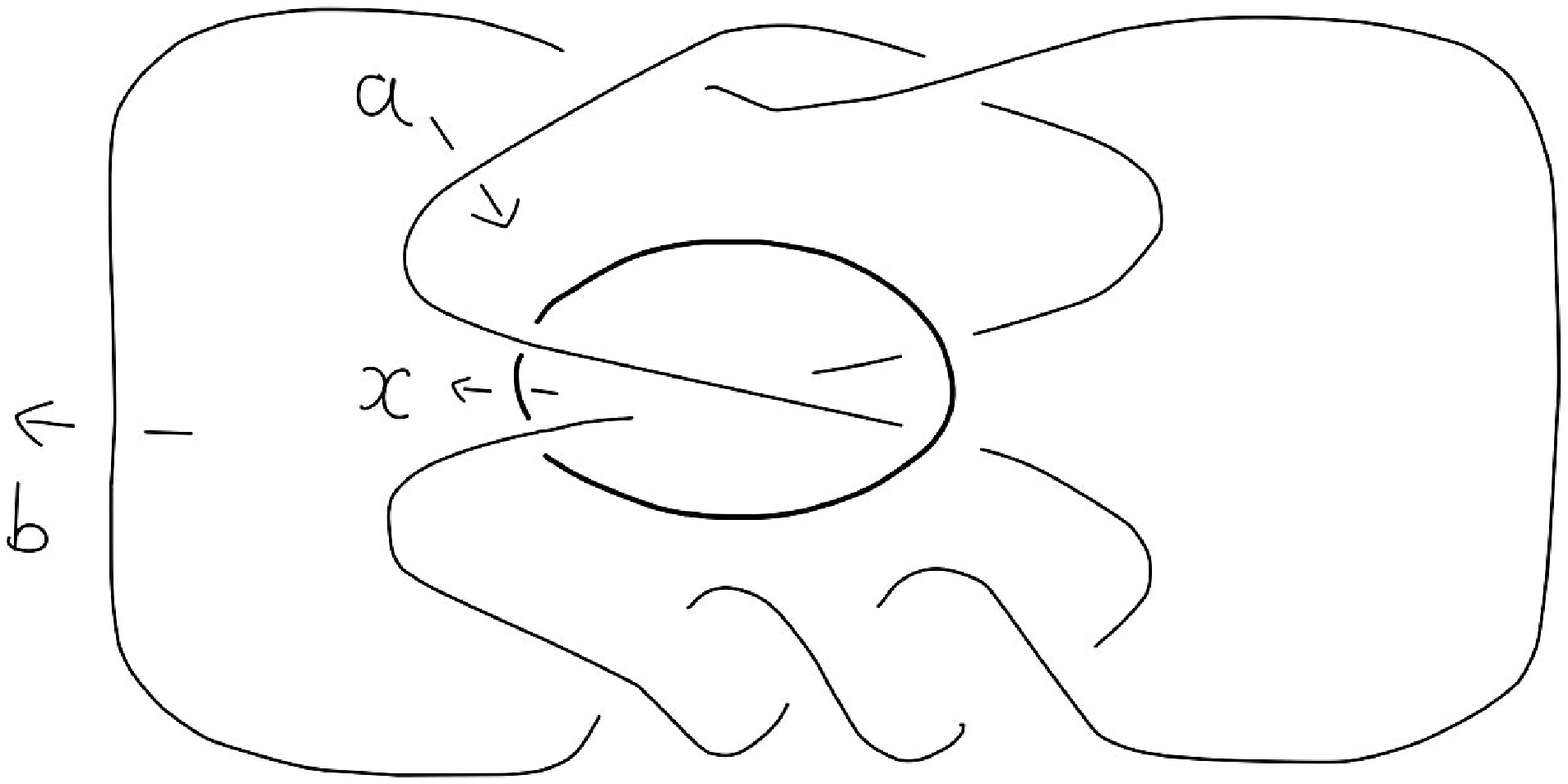}\caption{Link $L$}
  \end{center}
\end{figure}


\begin{prop}  \label{representation}
For a non-zero complex number $m$, there exists a representation $\rho_m : \pi_1(E_n) \to SL_2(\bC)$ such that
\begin{eqnarray*}
\rho_m(a)=
\left(
\begin{array}{cc}
 m & \displaystyle -\frac{\left(m^2-s\right) \left(s^{2 n+1}+1\right)}{m (s+1)} \\
 0 & \displaystyle m^{-1} \\
\end{array}
\right), \ \ 
\rho_m(b)=\displaystyle \frac{1}{s \alpha }
\left(
\begin{array}{cc}
 \beta &  \displaystyle -\frac{(s \alpha -m \beta)(m s \alpha -\beta )}{m  \beta } \\
 \beta &  \displaystyle \frac{m(m s \alpha -\beta ) +s \alpha }{m} \\
\end{array}
\right),
 \end{eqnarray*}
 and
\begin{eqnarray*}
\rho_m(x)=
 \left(
\begin{array}{cc}
 s^n & 0 \\
 \displaystyle \frac{s^n-s^{-n}}{s^{2 n+1}+1} & s^{-n} \\
\end{array}
\right),
\end{eqnarray*}
where $s$ is a solution to
 \begin{eqnarray}
0 =  m^8(s -  1\!\!\!&\!\!\!)(\!\!\!&\!\!\! s  + 1)^2 (s^{2 n}-s^2) s^{2 n+2}\\
-m^6\{s^{6 n+3}  \!\!\!&\!\!\!+ \!\!\!&\!\!\!(2 s^{6} + s^{5} - 4 s^{4} + s^{3} + s^{2} - s -1)s^{4 n+1} \nonumber \\
\!\!\!&\!\!\! - \!\!\!&\!\!\!(s^6 + s^{5} - s^{4} - s^{3} + 4 s^{2} - s -2)s^{2 n+2}   + s^{6}\} \nonumber\\
+m^4 \{(s^2 + \!\!\!&\!\!\!1 \!\!\!&\!\!\!)  s^{6 n+2} + (s^{6} + 2 s^{5} - 3 s^{4} - 2 s^{3} + 6 s^{2} - 4 s -2)s^{4 n+3}\nonumber\\
\!\!\!&\!\!\! - \!\!\!&\!\!\!  (2 s^{6} + 4 s^{5} - 6 s^{4} + 2 s^{3} + 3 s^{2} - 2 s -1)s^{2 n} + (s^2 + 1)s^5\} \nonumber\\
-m^2\{s^{6 n+3}  \!\!\!&\!\!\!+ \!\!\!&\!\!\!(2 s^{6} + s^{5} - 4 s^{4} + s^{3} + s^{2} - s -1)s^{4 n+1} \nonumber\\
\!\!\!&\!\!\! - \!\!\!&\!\!\!(s^6 + s^{5} - s^{4} - s^{3} + 4 s^{2} - s -2)s^{2 n+2}   + s^{6}\} \nonumber\\
+(s -1) (s \!\!\!&\!\!\!+ \!\!\!&\!\!\! 1  )^2 (s^{2 n}-s^2) s^{2 n+2} \nonumber
\end{eqnarray}
and $\alpha, \beta$ are given by
\begin{eqnarray*}
 \alpha \!\!&\!\! = \!\!& \!\! (s^2 - 1) s^{2 n} \{-m^6 (s - 1) s^2 (s^{2 n + 1} + 1) + m^4(s^{2 n+ 2} (s^4 - 2 s^2 + 3 s -1) + s^4 - 3 s^3 + 2 s^2 -1) \\
  && \!\! - m^2 s (s^ {2 n} (2 s^3 - s^2 + 1) - s (s^3 - s + 2)) + s^2 (s^{2 n} - s^2)\},\\
 \beta \!\!&\!\! =\!\! & \!\! m^7 s^{2 n + 2} (s^2 - 1)(s^3 + 1)  \\
   && \!\! - m^5 s^3 \{s^{4 n} (s^3 - s^2 + 1) + s^{2 n - 2} (s - 1) (s^3 + s + 1) (s^3 + s^2 + 1) - (s^3 - s + 1)\} \\
   && \!\! + m^3 s^2 (s^3 + 1) (s^{2 n} - 1) (s^{2 n}+ s^2)  - m s^3 (s^{2 n} - s^2) (s^{2 n} + s) .
\end{eqnarray*}
\end{prop}

In what follows, for simplicity, we denote the right hand side of $(1)$ by $r_0$.

\begin{proof}

For simplicity, put $A = \rho_m(a), \ B = \rho_m(b),\ X = \rho_m(x)$. By the aid of  Mathematica, we have
 \begin{eqnarray*}
 A X B A (X B)^{-1}  = 
\left(
\begin{array}{cc}
 s & 0 \\
 \displaystyle \frac{s^2-1}{s (s^{2n+1} + 1)} & \displaystyle \frac{1}{s} \\
\end{array}
\right)
+r_1
\left(
\begin{array}{cc}
 \displaystyle \frac{1}{m^3 s (s^{2n+1}+1) \alpha ^2} & \displaystyle -\frac{1}{m^3 s (s+1) \alpha ^2} \\
 \displaystyle \frac{s+1}{m^3 s^2 (s^{2n+1}+1)^2 \alpha ^2} & \displaystyle -\frac{1}{m^3 s^2 (s^{2n+1}+1) \alpha ^2} \\
\end{array}
\right),\\
  \end{eqnarray*} 
where
\begin{eqnarray*}
r_1 \!\!& \!\!= \!\!&\!\! -\alpha ^2 m s (m^2 s^{2n+2}-m^2-s^{2n+1}+s)
+\alpha  \beta  (m^2-1) (m^2+1) s^{2n+1} (s+1)\\
 &&\!\! +\beta ^2 m s^{2n} (m^2 s^{2n+1}-m^2 s-s^{2n+2}+1) \equiv 0 \mod r_0.
\end{eqnarray*}
Therefore, by $(1)$, we have $X=\{A X B A (X B)^{-1}\}^{n}$, that is, $\rho_m(x)=\rho_m \left(\{a x b a (x b)^{-1}\}^{n} \right)$. 

On the other hand, we can observe
\begin{eqnarray*}
  A X B \{A X B A (X B)^{-1} \} \equiv X B X^{-1}\{A X B A (X B)^{-1}\} X B   \mod r_0
\end{eqnarray*}
and so $A X B \{A X B A (X B)^{-1} \} = X B X^{-1}\{A X B A (X B)^{-1}\} X B$ by (1). Further more, we obtain
\begin{eqnarray*}
  X B \{A X B A (X B)^{-1} \}^{-1} (A X B)^{-1} X B 
  & = & X B (A X B \{A X B A (X B)^{-1} \})^{-1} X B \\
  & = & X B (X B X^{-1}\{A X B A (X B)^{-1}\} X B )^{-1} X B \\
  & = & \{A X B A (X B)^{-1}\}^{-1} X
\end{eqnarray*}
that is,
$\rho_m\left(\{a x b a (x b)^{-1}\}^{-1} x \right) = \rho_m\left(x b \{a x b a (x b)^{-1}\}^{-1} (a x b)^{-1} x b\right)$. 
This completes the proof.

\end{proof}

\begin{rem}
Since the representation $\rho_m$ comes from the holonomy representation obtained from the ideal triangulation of $E$ given in [TY],
the holonomy representation $\rho_m$ of $G(K_n)$ is given by the solution to $(1)$ which maximizes the hyperbolic volume of $S^3 \setminus K_n$.
\end{rem}

\section{Calculation of the twisted Alexander  polynomial}

The following is the main result of this paper.
\begin{thm}
The twisted Alexander polynomial of $K_n$ associated to $\rho_m$ is given by
\begin{eqnarray*}
\Delta_{K_n, \rho_m}(t) = 1 + \sum_{i=0}^{2n -1} \lambda_{i} (t^{i + 3} + t^{4n - i + 3} )+ t^{4 n + 6},
 \end{eqnarray*}
where 
 \begin{eqnarray*}
  \lambda_i = 
  \begin{cases}
  \displaystyle  \frac{(1 + m^2) (H s^{ i/2 +1 }\beta - s (s^{i/2 +1}- s^{-(i/2 +1)})  ( \eta_1 + \eta_2))}{H m \beta}  
  & {\rm if} \  0 \le i \le 2n - 2  \,\, {\rm and} \,\, i \,\, {\rm is \ even,}\\
    \displaystyle  \frac{s^{(i-1)/2} - s^{-(i-1)/2}}{s - s^{-1}} & {\rm if} \  0 \le i \le 2n - 2  \,\, {\rm and} \,\, i \,\, {\rm is \ odd,}\\
    \displaystyle  \frac{s^{n-1} - s^{-(n-1)}}{s - s^{-1}}  -\frac{(s^2 - 1) \eta_1}{H s^n \beta} & {\rm if} \  i = 2 n - 1
  \end{cases}  
 \end{eqnarray*}
and we put 
 \begin{eqnarray*}
 H \!\! & \!\!= \!\! & \!\! 1 - m^2 s + m^2 s^{2n +1} - s^{2n + 2},\\
  \eta_1\!\! & \!\!= \!\! & \!\! m \alpha - m s^{2n +1} \alpha + s^{2 n} \beta + m^2 s^{2 n} \beta,\\
  \eta_2\!\! & \!\!= \!\! & \!\! -m s \alpha + m s^{2n +1} \alpha - s^{2 n} \beta - s^{2n +1} \beta.
 \end{eqnarray*}
\end{thm}

To prove Theorem 3.1, it suffices to show

\begin{prop}  \label{mein prop} 
For simplicity, we put $S = s^n$ and $T = t^n$. The twisted Alexander polynomial $\Delta_{K_n, \rho_m}(t)$ is given by
\begin{eqnarray*}
&&\!\!\!\!\!\!\!\!\!\! \frac{S-T^2}{s-t^2}\frac{s}{S} \left( \frac{ m s - m S T^2 + (1 + m^2) (1 - s^2) S t T^2}{m (1 - s^2) t^2} +\frac{(1 + m^2) (1 - s S t^2 T^2) (\eta_1 + \eta_2)}{H m t^3 \beta} \right)\\
&+& \!\!\!\! \frac{1 - S T^2}{1 - s t^2} \frac{s}{S} \left( \frac{(1 + m^2) (1 - s^2) S - m S t + m s t T^2}{m (1 - s^2) t^3} -\frac{(1 + m^2) (s S - t^2 T^2) (\eta_1 + \eta_2)}{H m t^3  \beta} \right)\\
&+& \!\!\!\! \frac{1}{t^6} + T^4 + \frac{(1 - s^2) (1 + t^2) T^2 \eta_1}{H S t^4 \beta} .
\end{eqnarray*}
\end{prop}

By multiplying $t^6$ and rearranging with respect to $t$ , we obtain the formula of Theorem 3.1, when we use
\begin{eqnarray*}
\frac{S - T^2}{s - t^2} =\frac{S}{s} \sum_{i=0}^{n-1} \left( \frac{t^2}{s} \right)^i , \,\ 
\frac{S T^2 - 1}{s t^2 - 1} = \sum_{i=0}^{n-1} (s t^2)^i.
\end{eqnarray*}

\section{Proof of Proposition 3.2}
Recall that
 \begin{eqnarray*}
\pi_1 (E_n) \!\!&\!\! = \!\!&\!\! \langle a,b,x \ \vline \ \{ a x b a (x b)^{-1} \}^{-1} x = x b \{a x b a (x b)^{-1}\}^{-1} (a x b)^{-1} x b, \ 
x =  \{a x b a (x b)^{-1} \}^{n} \rangle\\
 \!\!& \!\!=\!\! &\!\! \langle a,c \ \vline \ ( a c a c^{-1} )^{n-1} = c (a c a c^{-1})^{-1} (a c)^{-1} c\rangle.
 \end{eqnarray*}

Then the twisted Alexander polynomial of $K_n$ is given by

\begin{eqnarray*}
\Delta_{K_n,\rho_m}(t) & = &  \frac{\displaystyle \det \Phi \left( \frac{\partial}{\partial a}( a c a c^{-1} )^{n-1}  -\frac{\partial}{\partial a} c (a c a c^{-1})^{-1} (a c)^{-1} c \right)}{\det \Phi(c-1)},
\end{eqnarray*}
where
\begin{eqnarray}
&& \Phi \left( \frac{\partial}{\partial a}( a c a c^{-1} )^{n-1}  -\frac{\partial}{\partial a} c (a c a c^{-1})^{-1} (a c)^{-1} c \right) \nonumber \\
&& = \sum_{i=1}^{n-1} t^{2(i-1)}  \rho_m \left( \left\{a x b a (x b)^{-1} \right\}^{i-1} \right) \left\{ \rho_m(1)+t^{2(n+1)} \rho_m(a x b)\right\} + t^{4n+1} \rho_m(x b x b a^{-1})  \\
&& + t^{2n-1} \rho_m \left(x b \left\{a x b a (x b)^{-1} \right\}^{-1}  \right) +t^{-3} \rho_m \left( x b \{a x b a (x b)^{-1}\}(a x b)^{-1} \right). \nonumber
\end{eqnarray}

For simplicity, we put
\begin{eqnarray*}
 \gamma_1 = s \alpha - m \beta \ ,\ 
 \gamma_2 = m s \alpha - \beta \ ,\ 
 \gamma_3 = m^2 s (s S^2 + 1) \alpha.
 \end{eqnarray*}
By the aid of Mathematica, the first term of the right hand side of $(2)$ is given by
\begin{eqnarray*}
&&\sum_{i=1}^{n-1} t^{2(i-1)} (AXBA(XB)^{-1})^{i-1} (E+t^{2(n+1)} AXB)\\
 &&   = 
\left(
\begin{array}{cc}
  \displaystyle \frac{(S T^2 - s t^2) (S t^2 \beta  T^2+m \alpha)}{m s t^2 (s t^2 - 1) \alpha } 
  & -\displaystyle \frac{T^2 (S T^2 - s t^2) (\gamma_1 \eta_2 + (m \alpha - \beta) \gamma_3)}{m^2 s (s+1) S \left(s t^2-1\right) \alpha  \beta } \\
  \displaystyle \frac{m C_1 \alpha -S t^2 T^2 C_2 \beta}{m s S (s S^2 + 1) t^2 (s - t^2)(s t^2 - 1) \alpha } 
  & \displaystyle \frac{C_3 t^4 T^4 + C_4 t^2 T^4 + C_5 t^6 T^2 + C_6 t^4 T^2 + C_7}{(s + 1) S^2 t^2 (s - t^2) (s t^2 -1)\gamma_3 \beta } \\
\end{array}
\right),
\end{eqnarray*}
where
\begin{eqnarray*}
 C_{1} \!\!&\!\!=\!\!&\!\! -t^4 s(s^2 - 1)S - T^2 \{t^2 (S^2 - s^4) - s(S^2 - s^2)\},\\
 C_{2} \!\!&\!\!=\!\!&\!\! -t^2(t^2 - 1) s(s + 1)S + T^2 \{t^2 (S^2 + s^3) + s(S^2 - s)\},\\
 C_{3} \!\!&\!\!=\!\!&\!\! (s^3 + S^2) \gamma_1 \eta_2 - \{s^3 (m s \alpha + \beta) - S^2 (m \alpha - \beta) \} \gamma_3,\\
 C_{4} \!\!&\!\!=\!\!&\!\!  - s (s + S^2) \gamma_1 \eta_2 + s \{s (m s \alpha + \beta) - S^2 (m \alpha - \beta)\} \gamma_3,\\
 C_{5} \!\!&\!\!=\!\!&\!\! - s (s + 1) S\{\gamma_1 \eta_2 + (\eta_1 + \eta_2 - (1 + m^2 S^2 - s S^2) \beta) \gamma_3 \},\\
 C_{6} \!\!&\!\!=\!\!&\!\! s (s + 1) S \{s \alpha \eta_2 - m (s + 1) S^2 \beta \gamma_2\},\\
 C_{7} \!\!&\!\!=\!\!&\!\! s (s + 1) S (s t^2 -1) (S t^2 - s T^2) \beta \gamma_3.
\end{eqnarray*}
Similarly, the second term of the right hand side of $(2)$ is given by
\begin{eqnarray*}
&& X B X B A^{-1}=
\left(
\begin{array}{cc}
 \displaystyle \frac{S^2 D_{1}}{\gamma_3 \alpha} 
 
 & \displaystyle \frac{m s D_1 D_2 - (s S^2 + 1) (s S^2 D_1 + m \gamma_3 \alpha) \beta^2}{ (s+1)\gamma_3 \alpha \beta ^2} \\
 
 \displaystyle \frac{(s+1) D_{2}}{ (s S^2+1)\gamma_3 \alpha} 
 
 & \displaystyle \frac{m s S^2 D_1 D_2 + s(s S^2 + 1) (m^2 s \alpha^2 - S^2 \beta^2) D_2}{S^2 (s S^2+1) \gamma_3 \alpha \beta ^2} -m \\
\end{array}
\right),
\end{eqnarray*}
where
\begin{eqnarray*}
D_1 \!\!&\!\!=\!\!&\!\! -(s + 1) \alpha \gamma_2 + m (\eta_1 + \gamma_2 + m S^2 \gamma_1)\beta,\\
D_2 \!\!&\!\!=\!\!&\!\! -\alpha \eta_2 + m S^2 (\eta_1 + m S^2 \gamma_1 + \gamma_2) \beta,
\end{eqnarray*}
the third term of the right hand side of (2) is given by
\begin{eqnarray*}
 X B \left\{AXBA(XB)^{-1}\right\}^{-1}=
\left(
\begin{array}{cc}
 \displaystyle \frac{S E_{1}}{m s \left(s S^2+1\right) \alpha  \beta } 
  &  \displaystyle -\frac{S \gamma_1 \gamma_2}{m \alpha  \beta } \\
  \displaystyle  \frac{(s+1) E_{2}}{m s S \left(s S^2+1\right)^2 \alpha  \beta } 
  &  \displaystyle \frac{E_{3}}{m S \left(s S^2+1\right) \alpha  \beta } \\ 
\end{array}
\right),
\end{eqnarray*}
where
\begin{eqnarray*}
 E_{1} \!\!&\!\!=\!\!&\!\! (s^2 - 1) \alpha \gamma_2 + m (\eta_1 + m S^2 \gamma_1 - s \gamma_2) \beta,\\
 E_{2} \!\!&\!\!=\!\!&\!\! (s - 1) \alpha \eta_2 + m S^2 (\eta_1 + m S^2 \gamma_1 -s \gamma_2) \beta,\\ 
 E_{3} \!\!&\!\!=\!\!&\!\! - s \alpha \eta_2 + m (s + 1) S^2 \beta \gamma_2,
\end{eqnarray*}
and the fourth term of the right hand side of (2) is given by
\begin{eqnarray*}
X B (A X B AXBA(XB)^{-1})^{-1}
 = 
\left(
\begin{array}{cc}
 \displaystyle  \frac{m F_{3} }{ \gamma_3 ^2 \beta ^2}  
 & \displaystyle  \frac{F_4}{m (s+1) \gamma_3 \alpha \beta ^2} \\ 
  \displaystyle \frac{m(s^2 - 1) F_{1} F_{2}}{S^2 (s  S^2+1) \gamma_3 ^2 \beta ^2}   
 & \displaystyle \frac{m F_5}{S^2 \gamma_3 ^2 \beta ^2} \\
\end{array}
\right),
\end{eqnarray*}
where
\begin{eqnarray*}
F_{1} \!\!&\!\!=\!\!&\!\! m (s + 1) S^2 (\eta_1 + m S^2 \gamma_1) \beta - \eta_2 \alpha ,\\
F_{2} \!\!&\!\!=\!\!&\!\! m (s + 1) S^2 (s S^2 + 1) \beta^2 - s F_1,\\
F_{3} \!\!&\!\!=\!\!&\!\! - \{m \beta (\eta_1 + m S^2 \gamma_1) + s \gamma_1 \gamma_2 - \gamma_2 \alpha\} F_2 + m s (s + 1) S^2 (s S^2 + 1) \gamma_1 \gamma_2 \beta^2,\\
 F_{4} \!\!&\!\!=\!\!&\!\! (s^2 - 1) \{m (\eta_1+m S^2 \gamma_1) \beta - \gamma_2 \alpha\} F_2 \\
&& + \gamma_3 \{m \gamma_2 \alpha - (m^2 \eta_1 + s^2 \eta_2 + m^3 S^2 \gamma_1 - s^2 (S^2 - 1) \gamma_2) \beta - 
 m s \gamma_1 \gamma_2\}\alpha,\\
 F_{5} \!\!&\!\!=\!\!&\!\! (s - 1) (s F_1 - m \gamma_3 \alpha) F_2 - m^2 S^2 (s S^2 + 1) \gamma_3 \alpha \beta^2.
\end{eqnarray*}
Therefore, the determinant of the right hand side of (2) is written as
\[
\frac{\sum_{i,j} U_{i,j} t^i T^j}{m^3 S^2 t^6 (s - t^2) (s t^2 - 1) \beta^2 \iota},
\]
where
\begin{eqnarray*}
  U_{0,0} \!\! & \!\!= \!\! & \!\! U_{4,0} = U_{6,0} = U_{2,4} = U_{10,4} = U_{6,8} = U_{8,8} = U_{12,8} 
   = -m^3 s S^2 \beta^2 \iota,\\ 
  U_{2,0}\!\! & \!\!= \!\! & \!\! U_{10,8} =  m^3 (s^2 + 1) S^2 \beta^2 \iota,\\
  H U_{3,0}\!\! & \!\! \equiv \!\! & \!\! H U_{9,8}  \equiv 
   -m^2 (m^2 + 1) s S^2 \beta (H s \beta - (s^2 - 1) (\eta_1 + \eta_2)) \iota  \mod r_0,\\  
  U_{5,0}\!\! & \!\! \equiv \!\! & \!\! U_{7,8}  \equiv m^2 (m^2 + 1) s S^2 \beta^2 \iota  \mod r_0,\\  
  H U_{1,2}\!\! & \!\! \equiv \!\! & \!\! H U_{11,6} \equiv  
  m^2 (m^2 + 1) (s - 1) s S \beta \eta_2 \iota  \mod r_0,\\  
  H U_{2,2}\!\! & \!\! = \!\! & \!\! H U_{6,2} = H U_{8,2} = H U_{4,6} \equiv H U_{6,6} = H U_{10,6} \equiv
  m^3 (s^2 - 1) s  S \beta \eta_1 \iota    \mod r_1,\\  
  H U_{3,2}\!\! & \!\! \equiv \!\! & \!\! H U_{9,6} \equiv  
  m^2 (m^2 + 1) (s - 1) S \beta \{H s S^2 \beta - s (s S^2 + 1) \eta_1 - (s^2 S^2 + s^2 + 1 ) \eta_2 \} \iota  \mod r_0,\\  
  H^2 U_{4,2}\!\! & \!\! \equiv \!\! & \!\! H^2 U_{8,6} \\
  \!\! & \!\! \equiv  \!\! & \!\!
  m (s - 1) s S  \{H^2 m^3 \alpha \beta + H (m^2 + 1) (m^2 s + s + 1) \beta \eta_2 - (m^2 + 1)^2 (s^2 - 1) \eta_2 (\eta_1 + \eta_2)\} \iota\\
  &&   \mod r_0,\\  
  H U_{5,2}\!\! & \!\! \equiv \!\! & \!\! H U_{7,6} \equiv  -m^2 (m^2 + 1) (s - 1) s S \beta \eta_2 \iota  \mod r_0,\\  
  H U_{7,2}\!\! & \!\! \equiv \!\! & \!\! H U_{5,6} \equiv  
  m^2 (m^2 + 1) (s - 1) s S \beta (H S^2 \beta - (s S^2 + 1) \eta_1 - (s S^2 - 1) \eta_2) \iota   \mod r_1,\\  
  H^2 U_{3,4}\!\! & \!\! \equiv \!\! & \!\! H^2 U_{9,4} \equiv  
  -m^2 (m^2 + 1) (s - 1)^2 s (s + 1) \eta_1 \eta_2 \iota    \mod r_0,\\  
  H^2 U_{4,4}\!\! & \!\! = \!\! & \!\! H^2 U_{8,4}\\
  \!\! & \!\! \equiv \!\! & \!\!
  m \{H^2 m^2 (s^2 - s + 1) S^2 \beta^2 + (m^2 + 1)^2 (s - 1)^2 s \eta_2 (-H S^2 \beta + (s S^2 + 1) \eta_1 + s S^2 \eta_2)\} \iota \\
  &&     \mod r_1,\\  
  H^2 U_{5,4}\!\! & \!\! \equiv\!\! & \!\! H^2 U_{7,4}\\  
  \!\! & \!\! \equiv \!\! & \!\! 
  -(m^2 + 1) (s - 1) s \{(s - 1) \eta_2 (m^3 H \alpha + (m^2 + 1) \eta_2) + m^2 S^2 H \beta (H \beta - (s + 1) (\eta_1 + \eta_2))\} \iota  \\
  &&         \mod r_0,\\  
  H^2 U_{6,4}\!\! & \!\! \equiv \!\! & \!\! 
  -2 m s (H m S \beta - (m^2 + 1) (s - 1) \eta_2) (H m S \beta + (m^2 + 1) (s - 1) \eta_2) \iota  \mod r_0,
\end{eqnarray*}
where we put
$\iota = m^2 s^2 (s + 1) S (s S^2 + 1)^3 \alpha^3 \beta$, and the other $U_{i,j}$'s are $0$.

On the other hand, by the aid of Mathematica,
\begin{eqnarray*}
\det \Phi(c-1)  \!\!&\!\! = \!\!&\!\! 
\det \left( t^{2n+1}\rho_m(x b)-
\left(
\begin{array}{cc}
1 & 0\\
0 & 1
\end{array}
\right) \right)\\
\!\!&\!\! = \!\!&\!\!
\frac{m S H \beta + m S H t^2 T^4 \beta - (m^2 + 1) (s - 1) t T^2 \eta_2}{m S H \beta} 
-\frac{(S^2-1) t T^2}{m S (s S^2+1) H \alpha \beta}r_1\\
\!\!&\!\! = \!\!&\!\!
\frac{m S H \beta + m S H t^2 T^4 \beta - (m^2 + 1) (s - 1) t T^2 \eta_2}{m S H \beta} .
\end{eqnarray*}

Consequently, we have
\begin{eqnarray}
 \Delta_{K_n, \rho_m}(t) = \frac{\sum_{i,j} V_{i,j} t^i T^j}{H m^2 S t^6 (s - t^2) (s t^2 - 1) \beta},
\end{eqnarray}
where
\begin{eqnarray*}
V_{0, 0} \!\! & \!\!= \!\! & \!\! V_{4, 0} = V_{6, 0} = V_{4, 4} = V_ {6, 4} = V_{10, 4} = -H m^2 s S \beta,\\
V_{2, 0} \!\! & \!\!= \!\! & \!\! V_{8, 4} = H m^2 (s^2 + 1) S \beta,\\
V_{3, 0} \!\! & \!\!= \!\! & \!\! V_{7, 4} = m (m^2 + 1) s S \{(s^2 - 1)( \eta_1 + \eta_2)  - H s \beta\},\\
V_{5, 0} \!\! & \!\!= \!\! & \!\! V_{5, 4} = H m (m^2 + 1) s S \beta,\\
V_{2, 2} \!\! & \!\!= \!\! & \!\! V_{8, 2} = m^2 s (s^2 - 1) \eta_1,\\
V_{3, 2} \!\! & \!\!= \!\! & \!\! V_{7, 2} = m (m^2 + 1) (s -1) s \{(s + 1) \eta_1 + \eta_2 \}\\
V_{4, 2} \!\! & \!\!= \!\! & \!\! V_{6, 2} = (s - 1) s \{(m^2 + 1) \eta_2 + H m^3 \alpha\},\\
V_{5, 2} \!\! & \!\!= \!\! & \!\! -2 m (m^2 + 1) (s - 1) s \eta_2,
\end{eqnarray*}
and the other $V_{i,j}$'s are $0$.
By the aid of Mathematica, the difference between the right hand side of (3) and the formula in Proposition 3.2 is equal to 
\begin{eqnarray*}
\frac{s \zeta_1 + t \zeta_2 - 2 t^2 \zeta_1 + t^3 \zeta_2 + s t^4 \zeta_1}{H m^2 S t^3 (s + 1) (s - t^2) (s t^2 - 1) \beta}  T^2,
\end{eqnarray*} 
where
\begin{eqnarray*}
\zeta_1 \!\! & \!\!= \!\! & \!\! m (m^2 + 1) s (s + 1) (H S^2 \beta - s (S^2 - 1) \eta_1 - (s S^2 - 1) \eta_2), \\
\zeta_2 \!\! & \!\!= \!\! & \!\! H m^2 s (m \alpha - m s^2 \alpha + s \beta + S^2 \beta) - (s^2 - 1) (m^2 \eta_1 + m^2 s^3 \eta_1 + s \eta_2 + m^2 s \eta_2).
\end{eqnarray*} 
Note that $\zeta_1 = 0$ by the definition of $H, \eta_1$ and $\eta_2$ and that
 \begin{eqnarray*}
\zeta_2 =
m \{(m^2 (s^2 - s + 1) - s) (s^3 S^2 + 1) - H s (s - 1)\} r_0 = 0.
\end{eqnarray*} 
This completes the proof of Proposition 3.2.

\vspace{4mm}

\end{document}